\documentclass[a4paper,12pt]{article}
\usepackage{amssymb,amsmath}

\def\sectsign{\mathhexbox278}
\newcommand{\GL}{\mathrm{GL}}
\renewcommand\epsilon{\varepsilon}
\def\Q{\mathbb Q}
\def\Z{\mathbb Z}

\def\rad{\mathop{\rm rad}}
\def\Ext{\mathop{\rm Ext}\nolimits^1}
\def\Hom{\mathop{\rm Hom}\nolimits^{\vphantom{1}}}
\def\diag{\mathop{\rm diag}\nolimits}

\def\led{\vartriangleleft}
\def\ged{\vartriangleright}
\def\ledeq{\trianglelefteq}
\def\gedeq{\trianglerighteq}

\def\notged{\ntriangleright}

\def\proof{{\bf Proof. }}
\def\endproof{\hfill$\square$}

\renewcommand{\le}{\leqslant}
\renewcommand{\ge}{\geqslant}

\makeatletter
\def\@begintheorem#1#2{\it \trivlist
       \item[\hskip \labelsep{\bf #1\ #2.}]}
\def\@opargbegintheorem#1#2#3{\it \trivlist
       \item[\hskip \labelsep{\bf #1\ #2.\ (#3)}]}
\makeatother

\newtheorem{theorem}{Theorem}[section]
\newtheorem{lemma}[theorem]{Lemma}

\newtheorem{definition}[theorem]{Definition}
\newtheorem{proposition}[theorem]{Proposition}

\begin{document}

\begin{center}
 \Large On some extensions of $p$-restricted completely splittable $\GL(n)$-modules
\normalsize
\vskip 1.5em%
Vladimir Shchigolev\\
{\small Lomonosov Moscow State University,\\ e-mail: shchigolev\_vladimir@yahoo.com}
\end{center}

\begin{abstract}
In this paper, we calculate the space
$\Ext_{\GL(n)}(L_n(\lambda),L_n(\mu))$, where $\GL(n)$ is the general
linear group of degree $n$ over an algebraically closed field of
positive characteristic, $L_n(\lambda)$ and $L_n(\mu)$ are rational
irreducible $\GL(n)$-modules
with highest weights $\lambda$ and $\mu$ respectively,
the restriction of $L_n(\lambda)$ to any Levi subgroup of $\GL(n)$
is semisimple, $\lambda$ is a $p$-restricted weight and
$\mu$ does not strictly dominate $\lambda$.
\end{abstract}

\section{Introduction}

In this paper, we fix an algebraically closed field $K$
of characteristic $p>0$.
We denote by $\Sigma_n$ the symmetric group of order $n$
interpreted in this paper as the group of all bijections of
the set $\{1,\ldots,n\}$. Let $\GL(n)$ denote the group of
all invertible $n\times n$-matrices over $K$.

A {\it partition} of an integer $n$ is an infinite countable
weakly decreasing sequence of nonnegative integers such that
the sum of its elements equals $n$.
The {\it height} of a partition $\lambda$ is the number $h(\lambda)$ of
its nonzero entries.
A partition $\lambda$ that does not contain $p$ or more identical entries is
called {\it $p$-regular}.
A partition $\lambda$ is said to dominate a partition $\mu$
if $\sum_{j=1}^i\lambda_j\ge\sum_{j=1}^i\mu_j$ for any $i\ge1$.
This fact is denoted by $\lambda\gedeq\mu$.
The formula $\lambda\ged\mu$ means that $\lambda\gedeq\mu$ and $\lambda\ne\mu$.
Let $\lambda^t$ denote the partition conjugate to $\lambda$.
That is the Young diagram of $\lambda^t$ is the transposed
Young diagram of $\lambda$.

To each partition $\lambda$ of $n$, there corresponds a $K\Sigma_n$-module
$S^\lambda$, which is called the {\it Specht module}
(see, for example, \cite[Definition~4.3]{James1}).
Let $D^\lambda=S^\lambda/\rad S^\lambda$ for a $p$-regular partition $\lambda$.
The map $\lambda\mapsto D^\lambda$ defines a one-to-one correspondence between
$p$-regular partitions of $n$ and irreducible $K\Sigma_n$-modules.

By a {\it Young subgroup} of $\Sigma_n$ we mean any subgroup
$$
\{\sigma\in\Sigma_n:\sigma A_i=A_i,\;i=1,\dots,k\},
$$
where $A_1\cup\cdots\cup A_k=\{1,\ldots,n\}$ is an arbitrary
disjoint partition.

\begin{definition}[\mbox{\cite[Definition~0.1]{Kleshchev_notecs}}]
An irreducible $K\Sigma_n$-module
is called {\it completely splittable\/}
if and only if its restriction to any Young subgroup of $\Sigma_n$
is semisimple.
\end{definition}

For any partition $\lambda$, let $\chi(\lambda)=\lambda_1-\lambda_{h(\lambda)}+h(\lambda)$
if $\lambda$ is nonzero and $\chi(\lambda)=0$ otherwise.
A criterion for complete splittability of $D^\lambda$ in terms of $\lambda$
was obtained by A.S.\,Kleshchev in~\cite{Kleshchev_notecs}.

\begin{proposition}[\mbox{\cite[Theorem~2.1]{Kleshchev_notecs}}]\label{proposition:0.125}
$D^\lambda$ is completely splittable if and only if $\chi(\lambda)\le p$.
\end{proposition}

In what follows a partition $\lambda$ such that
$D^\lambda$ is completely splittable is also called
{\it completely splittable}.
A completely splittable partition containing more than one nonzero parts
and at least one rim $p$-hook is called {\it big}
(see~\cite[Definition~4]{Shchigolev11}).

The result for the general linear group
similar to Proposition~\ref{proposition:0.125}
was proved in~\cite{Kleshchev_gjs11}.
In this connection, we are interested to
prove an analog of the following result for the general linear group.

\begin{proposition}[\mbox{\cite[Theorem~6 and Corollary~6]{Shchigolev11}}]\label{proposition:0.25}
Let $p>2$, $\lambda$ be a completely splittable partition and
$\mu$ be a $p$-regular partition of $n$ such that $\lambda\notged\mu$.
Then
$$
\Ext_{\Sigma_n}(D^\lambda,D^\mu)\cong
\left\{
\begin{array}{ll}
K &\mbox{if }\lambda\mbox{ is big and }\mu=\tilde\lambda;\\
0 &\mbox{otherwise}.
\end{array}
\right.
$$
If $\lambda$ is big then $\rad S^\lambda$ is a nonzero homomorphic
image of $S^{\tilde\lambda}$. Otherwise $\rad S^\lambda=0$.
\end{proposition}

The main aim of the present paper is to prove the above mentioned
analog for $p$-restricted completely splittable weights
(see the definitions of~\sectsign 2).
In fact, the desired result (Theorem~\ref{theorem:1})
follows from Proposition~\ref{proposition:0.25} by applying
the Schur functor and the restriction functors $R^{n+m}_n$
from~\cite{Kleshchev_lr10}.

Note that the question remains open for not $p$-restricted
completely splittable weights.

\section{Preliminaries}

For an integer $n\ge0$, let $X(n)$ denote the
set of all $n$-tuples of integers.
The subset of $X(n)$ consisting of sequences $(\lambda_1$, \ldots, $\lambda_n)$
such that $\lambda_1\ge\cdots\ge\lambda_n$ is denoted by $X^+(n)$.
Elements of $X(n)$ are called {\it weights} and
elements of $X^+(n)$ are called {\it dominant weights}.
Let $\Lambda^+(n,m)$ denote the subset of $X^+(n)$
consisting of all sequences $(\lambda_1$, \ldots, $\lambda_n)$ such that
$\lambda_1\ge\cdots\ge\lambda_n\ge0$ and $\lambda_1+\cdots+\lambda_n=m$.

A vector $v$ of a rational $\GL(n)$-module $M$ is said
to have weight $\lambda\in X(n)$
if $\diag(t_1,\ldots,t_n)v=t_1^{\lambda_1}\cdots t_n^{\lambda_n}v$,
for all $t_1$, \ldots, $t_n\in K\setminus\{0\}$.
We denote by $M^\lambda$ the subspace of $M$
consisting of all vectors of weight $\lambda$.
We say that $\lambda$ is a weight of $M$ if $M^\lambda\ne0$.

We fix the weights $\epsilon_i=(0^{i-1},1,0^{n-i})$ and
$\alpha_{i,j}=\epsilon_i-\epsilon_j$, where $1\le i,j\le n$.
We introduce the following partial dominance order on weights:
$\mu\le\lambda$ if $\lambda-\mu$ is a sum of $\alpha_{i,j}$ for $i<j$
(possibly empty and with repetitions).
As usual, $\mu<\lambda$ means that $\mu\le\lambda$ and $\mu\ne\lambda$.
A dominant weight $\lambda\in X^+(n)$ is called
$p$-restricted if $\lambda_i-\lambda_{i+1}<p$ for $1\le i<n$.

For a dominant weight $\lambda\in X^+(n)$, row $i$ is called {\it removable}
if $1\le i<n$ and $\lambda_i>\lambda_{i+1}$.
If a weight $\lambda\in X^+(n)$ has at least one removable row,
then we put $\psi_n(\lambda)=j-i+\lambda_i-\lambda_{j+1}$,
where $i$ and $j$ are the smallest and the largest removable rows
of $\lambda$ respectively.
If $\lambda$ does not have removable rows, i.e. all entries of $\lambda$
are identical, then we put $\psi_n(\lambda)=0$.

Let $[V:D]$ denote the composition multiplicity of an irreducible module $D$
in a module $V$, when it is clear over what ring both
modules are considered.

For $\lambda\in X^+(n)$, let $\Delta_n(\lambda)$ and $L_n(\lambda)$ denote
the {\it Weyl module} and the {\it irreducible module}
with highest weight $\lambda$ respectively
(see the definitions from~\cite{Jantzen1}).
As we distinguish between partitions and
weights, the following notation is helpful.
Let $a$ be a sequence of length $n$ and $m$ be such that
$a_i\ne0$ implies $i\le m$.
Then we denote by $(a)_m$ the sequence
$(a_1$, \ldots, $a_m)$ if $m<n$ and the sequence
$(a_1$, \ldots, $a_n,0^{m-n})$ if $m\ge n$.
In this definition $n$ and $m$ are nonnegative integers or $+\infty$.
We shall mean the following abbreviations:
$$
S^\lambda=S^{(\lambda)_{+\infty}},         \quad
\lambda^t=((\lambda)_{+\infty})^t,         \quad
\Delta_n(\lambda)=\Delta_n((\lambda)_n),   \quad
L_n(\lambda)=L_n((\lambda)_n)
$$
always, when the corresponding modules and sequences are well defined.

A {\it standard Levi subgroup} of $\GL(n)$ is any subgroup of the form
$\GL(n_1)\times\cdots\times\GL(n_k)$, where $n_1+\cdots+n_k=n$.
A subgroup of $\GL(n)$ is called a {\it Levi subgroup} of $\GL(n)$
if it is conjugate to a standard Levi subgroup.
An irreducible module $L_n(\lambda)$ is called {\it completely splittable} if
its restriction to any Levi subgroup of $\GL(n)$ is semisimple.
A weight $\lambda\in X^+(n)$ such that the module $L_n(\lambda)$
is completely splittable is also called {\it completely splittable}.

In~\cite{Kleshchev_gjs11}, a criterion for complete splittability of $\lambda$
is proved, which in the case of $p$-restricted weights
takes the following form.
\begin{proposition}\label{proposition:0.375}
A $p$-restricted weight $\lambda\in X^+(n)$ is completely splittable
if and only if $\psi_n(\lambda)\le p$.
\end{proposition}

In the proof of the main result, we shall
use the following fact that follows from~\cite[II.2.14(4)]{Jantzen1}.
\begin{proposition}\label{proposition:0.5}
Let $\lambda,\mu\in X^+(n)$ and $\lambda\not<\mu$. Then
$$
\Ext_{\GL(n)}(L_n(\lambda),L_n(\mu))\cong\Hom_{\GL(n)}(\rad\Delta_n(\lambda),L_n(\mu)).
$$
\end{proposition}
Here and in what follows cohomologies of algebraic groups
are understood as rational cohomologies (see~\cite[I.4]{Jantzen1}).

In the proof of Theorem~\ref{theorem:1}, we use the hyperalgebra $U(n)$
for $\GL(n)$.
It has the advantage of being generated by elements taking
homogeneous vectors to homogeneous vectors.
Below is its exact definition.

Let ${\mathfrak gl}(n,\Q)$ be the Lie algebra of all $n\times n$-matrices
over the field of rational numbers $\Q$ with respect to commutation and
$U(n,\Q)$ be the universal enveloping algebra of ${\mathfrak gl}(n,\Q)$.
Let $X_{i,j}$ denote the matrix having $1$ at the intersection of
row $i$ and column $j$ and zeros elsewhere.
Thus $X_{i,j}\in{\mathfrak gl}(n,\Q)\subset U(n,\Q)$
Let $U(n,\Z)$ denote the $\Z$-subalgebra of $U(n,\Q)$
generated by the elements
$$
\begin{array}{l}
X_{i,j}^{(r)}:=\frac{(X_{i,j})^r}{r!}\mbox{ for integers }1\le i,\; j\le m, i\ne j\mbox{ and }r\ge0;\\[12pt]
\binom{X_{i,i}}{r}:=\frac{X_{i,i}(X_{i,i}-1)\cdots(X_{i,i}-r+1)}{r!}\mbox{ for integers }1\le i\le m\mbox{ and }r\ge0.
\end{array}
$$
As is shown in~\cite{Jantzen1}, any finite dimensional rational
$\GL(n)$-module can be interpreted as a $U(n)$-module.
The subalgebras of $U(n)$ generated by all elements
$X_{i,j}^{(r)}$ for $i<j$ is denoted by $U^+(n)$.
The subalgebras of $U(n)$ generated by all elements
$X_{i,j}^{(r)}$ for $i>j$ is denoted by $U^-(n)$.
The subalgebra of $U(n)$ generated by all elements $\binom{X_{i,i}}{r}$
is denoted by $U^0(n)$.
The decomposition $U(n)=U^-(n)U^{0}(n)U^+(n)$ is proved in~\cite{Jantzen1}.
Finally, let us call a vector $v$ of some $U(n)$-module {\it primitive}
if $U^+(n)v=Kv$.

\section{Main result}

We shall prove the following simple fact.
\begin{proposition}\label{proposition:1}
Let $S$ be a $K$-algebra,
$V$ be an $S$-module, $L$ be an irreducible $S$-module and
$e\in S$ be an idempotent such that $eL\ne0$.
Then $\Hom_S(V,L)$ is isomorphic to a $K$-subspace of $\Hom_{eSe}(eV,eL)$
\end{proposition}
\proof
Consider the $K$-linear map from $\Hom_S(V,L)$ to $\Hom_{eSe}(eV,eL)$
defined by $\varphi\mapsto\varphi|_{eV}$.
We claim that it is an embedding.
Let $\varphi$ be a nonzero element of $\Hom_S(V,L)$.
By hypothesis, $ex\ne0$ for some $x\in L$.
Since $L$ is irreducible, we have $x=\varphi(v)$ for some $v\in V$.
Hence $\varphi(ev)=ex\ne0$ and $\varphi|_{eV}\ne0$.
\endproof

\begin{lemma}\label{lemma:1}
Let $\lambda,\mu\in\Lambda^+(n,m)$, $m\le n$, $\lambda\not<\mu$ and
$\lambda^t,\mu^t$ be $p$-regular.
Then
$\dim\Ext_{\GL(n)}(L_n(\lambda),L_n(\mu))\le\dim\Ext_{\Sigma_m}(D^{\lambda^t},D^{\mu^t})$.
\end{lemma}
\proof
It follows from Proposition~\ref{proposition:0.5} that
$\Ext_{\GL(n)}(L_n(\lambda),L_n(\mu))\cong\Hom_{S(n,m)}(\rad\Delta_n(\lambda),L_n(\mu))$,
where $S(n,m)$ is the Schur algebra (see~\cite[\sectsign2.3]{J.A.Green1}).
Proposition~\ref{proposition:1} applied to
$e:=\xi_{u,u}$, where $u=(1,2,\ldots,m,0^{n-m})$
(see~\cite[\sectsign\sectsign2.3,6]{J.A.Green1}),
$V:=\rad\Delta_n(\lambda)$, $L:=L_n(\mu)$ and $S:=S(n,m)$
yields
$$
\dim\Hom_{S(n,m)}(\rad\Delta_n(\lambda),L_n(\mu))\le
\dim\Hom_{\Sigma_m}(\rad S^{\lambda^t},D^{\mu^t}).
$$
Applying $\Hom_{\Sigma_m}(-,D^{\mu^t})$ to the exact sequence
$0\to\rad S^{\lambda^t}\to S^{\lambda^t}\to D^{\lambda^t}\to0$,
we get the exact sequence
$$
\begin{array}{c}
\hbox to 0.95\textwidth{$0\to\Hom_{\Sigma_m}(D^{\lambda^t},D^{\mu^t})
\to\Hom_{\Sigma_m}(S^{\lambda^t},D^{\mu^t})$\hfil}\\[4pt]
\hbox to 0.95\textwidth{\hfil$\to\Hom_{\Sigma_m}(\rad S^{\lambda^t},D^{\mu^t})
\stackrel\iota\to\Ext_{\Sigma_m}(D^{\lambda^t},D^{\mu^t}).$}
\end{array}
$$
Hence $\iota$ is an embedding and
$\dim\Hom_{\Sigma_m}(\rad S^{\lambda^t},D^{\mu^t})\le\dim\Ext_{\Sigma_m}(D^{\lambda^t},D^{\mu^t})$.
\endproof

\begin{definition}\label{definition:1}
Let $\lambda$ be a $p$-restricted weight of $X^+(n)$.
We call $\lambda$ big if $\lambda$ has a removable row and
\begin{enumerate}
\itemsep=0pt
\item $j-i+\lambda_i-\lambda_{j+1}\le p$ (i.e. $\lambda$ is completely splittable);
\item $\lambda_i-\lambda_{j+1}>1$,
      $j-1+\lambda_i-\lambda_{j+1}\ge p$,
      $i-\lambda_i+\lambda_{j+1}+p+1\le n$,
\end{enumerate}
where $i$ and $j$ are the smallest and the largest removable rows
of $\lambda$ respectively.
For a big $\lambda\in X^+(n)$, we put
$$\hat\lambda=\lambda-\sum_{k=1}^{i-j-\lambda_i+\lambda_{j+1}+p+1}\alpha_{j+\lambda_i-\lambda_{j+1}-p-1+k,j+k}.$$
\end{definition}
In other words, a $p$-restricted weight $\lambda$ is big if
$\lambda$ is completely splittable and one can nontrivially move
$i-j-\lambda_i+\lambda_{j+1}+p+1$ nodes from the last column to column $\lambda_{j+1}+1$
to obtain a weight of $X^+(n)$, which is denoted by $\hat\lambda$.

{\it Example.}
Let $p=7$ and $\lambda=(5,5,5,4,3,1,1,1)$.
Then $\lambda$ is big and $\hat\lambda=(5,4,4,4,3,2,2,1)$.

\begin{theorem}\label{theorem:1}
Let $\lambda$ be a completely splittable $p$-restricted weight of $X^+(n)$ and
$\mu\in X^+(n)$ such that $\lambda\not<\mu$.
Then
$$
\Ext_{\GL(n)}(L_n(\lambda),L_n(\mu))\cong
\left\{
\begin{array}{ll}
K &\mbox{if }\lambda\mbox{ is big and }\mu=\hat\lambda;\\
0 &\mbox{otherwise}.
\end{array}
\right.
$$
If $\lambda$ is big then $\rad\Delta_n(\lambda)$ is a nonzero homomorphic
image of $\Delta_n(\hat\lambda)$. Otherwise $\rad\Delta_n(\lambda)=0$.
\end{theorem}
\proof
Multiplying the modules $L_n(\lambda)$ and $L_n(\mu)$
by a suitable power of
the representation $\det_n=\Delta_n(1^n)$,
we may assume that $\lambda_n=0$ (see Proposition~\ref{proposition:0.5}).
We put $m=\sum_{i=1}^n\lambda_i$.

Case 1: $\lambda\not>\mu$.
By Proposition~\ref{proposition:0.5}
we get $\Ext_{\GL(n)}(L_n(\lambda),L_n(\mu))=0$.
On the other hand, for a big $\lambda$, we have $\lambda>\hat\lambda$.

Case 2: $\lambda>\mu$ and $n\ge m$.
If $p=2$ then we would have $\lambda=(1^m,0^{n-m})$ and
the conditions $\mu<\lambda$ and $\mu\in X^+(n)$ would not hold simultaneously.
Hence $p>2$.
We have $h(\lambda^t)<p$, $\lambda^t\led\mu^t$ and $h(\mu^t)<p$.
Therefore $\mu^t$ and $\lambda^t$ are $p$-regular.
By Lemma~\ref{lemma:1}, we get
\begin{equation}\label{equation:0.5}
\dim\Ext_{\GL(n)}(L_n(\lambda),L_n(\mu))\le\dim\Ext_{\Sigma_m}(D^{\lambda^t},D^{\mu^t}).
\end{equation}

In the case under consideration,
$\lambda^t$ is a big partition if and only if $\lambda$ is a big weight;
if these conditions are satisfied, then $\widetilde{\lambda^t}=\hat\lambda^t$.
It follows from these facts, inequality (\ref{equation:0.5})
and Proposition~\ref{proposition:0.25} that
\begin{equation}\label{equation:1}
\dim\Ext_{\GL(n)}(L_n(\lambda),L_n(\mu))\le
\left\{
\begin{array}{ll}
1&\mbox{ if }\lambda\mbox{ is big and }\mu=\hat\lambda;\\
0&\mbox{ otherwise}.
\end{array}
\right.
\end{equation}
Therefore it remains to prove that for a big $\lambda$,
the module $\rad\Delta_n(\lambda)$ is a nonzero homomorphic
image of $\Delta_n(\hat\lambda)$.

By~\cite[6.6b]{J.A.Green1} and Proposition~\ref{proposition:0.25}, we have
\begin{equation}\label{equation:2}
[\rad\Delta_n(\lambda):L_n(\hat\lambda)]=[\Delta_n(\lambda):L_n(\hat\lambda)]=
[S^{\lambda^t}:D^{\hat\lambda^t}]=[S^{\lambda^t}:D^{\widetilde{\lambda^t}}]=1.
\end{equation}
Therefore there is a nonzero vector $v\in\rad\Delta_n(\lambda)$
of weight $\hat\lambda$.
Let $V$ be a $U(n)$-submodule in $\rad\Delta_n(\lambda)$
generated by $v$, $\gamma\in X(n)^+$ be a maximal
with respect to the order $\ge$ weight of $V$
and $w$ be a nonzero element of $V^\gamma$.
Clearly, $w$ is a primitive vector and $w\in U^+(n)v$.
Hence $\Hom_{\GL(n)}(\Delta_n(\gamma),\Delta_n(\lambda))\ne0$ and
$\gamma\ge\hat\lambda$.
Therefore by Proposition~\ref{proposition:0.25}, we get
$$
0\ne[\rad\Delta_n(\lambda):L_n(\gamma)]=[\rad S^{\lambda^t}:D^{\gamma^t}]\le
[S^{\widetilde{\lambda^t}}:D^{\gamma^t}].
$$
Hence $\hat\lambda^t=\widetilde{\lambda^t}\ledeq\gamma^t$, $\hat\lambda\ge\gamma$ and
$\hat\lambda=\gamma$.
Now it follows from~(\ref{equation:1}) and~(\ref{equation:2})
that the image of an arbitrary nonzero homomorphism from
$\Delta_n(\hat\lambda)$ to $\Delta_n(\lambda)$ coincides with $\rad\Delta_n(\lambda)$.

Case 3: $\lambda>\mu$.
We put $\nu=(\lambda_1$, \ldots, $\lambda_n,0^m)$.
First suppose that $\lambda$ is big.
Then $\nu$ is also big, $\hat\nu=(\hat\lambda_1$, \ldots, $\hat\lambda_n,0^m)$,
and by case~2 there exists an exact sequence
$$
\Delta_{n+m}(\hat\nu)\to\Delta_{n+m}(\nu)\to L_{n+m}(\nu)\to0.
$$
Applying the functor $R^{n+m}_n$ defined in~\cite[\sectsign2]{Kleshchev_lr10}
to this sequence and taking into account~\cite[2.3]{Kleshchev_lr10},
we get an exact sequence
$$
\Delta_n(\hat\lambda)\to\Delta_n(\lambda)\to L_n(\lambda)\to0.
$$
Therefore $\rad\Delta_n(\lambda)$ is a homomorphic image of $\Delta_n(\hat\lambda)$.
In view of Proposition~\ref{proposition:0.5}, it remains to show that
$\rad\Delta_n(\lambda)$ is nonzero.
By case~2, we have $1=[\Delta_{n+m}(\nu):L_{n+m}(\hat\nu)]$.
Applying the functor $R^{n+m}_n$ and
taking into account~\cite[2.3]{Kleshchev_lr10},
we get $1=[\Delta_n(\lambda):L_n(\hat\lambda)]=[\rad\Delta_n(\lambda):L_n(\hat\lambda)]$.

Now suppose that $\lambda$ is not big.
In view of Proposition~\ref{proposition:0.5},
we only need to prove that $\rad\Delta_n(\lambda)=0$.

If $\nu$ is also not big, then $\Delta_{n+m}(\nu)\cong L_{n+m}(\nu)$
by case~2.
Applying the functor $R^{n+m}_n$ and
taking into account~\cite[2.3]{Kleshchev_lr10},
we get $\Delta_n(\lambda)\cong L_n(\lambda)$ and $\rad\Delta_n(\lambda)=0$.

If $\nu$ is big, then $n<m$, $\hat\nu_{n+1}>0$ and by case~2
there exists an exact sequence
$$
\Delta_{n+m}(\hat\nu)\to\Delta_{n+m}(\nu)\to L_{n+m}(\nu)\to0.
$$
Applying the functor $R^{n+m}_n$ to this sequence and taking
into account~\cite[2.4 and 2.3]{Kleshchev_lr10},
we get $\Delta_n(\lambda)\cong L_n(\lambda)$ and $\rad\Delta_n(\lambda)=0$.
\endproof


\begin{thebibliography}{1}\label{bibl}

\bibitem{Kleshchev_lr10}
J. Brundan and A.S. Kleshchev, Modular Littlewood-Richardson coefficients, {\em
  Math. Z.}, {\bf 232} (1999), 287--320.

\bibitem{Kleshchev_gjs11}
J. Brundan, A.S. Kleshchev, and I.D. Suprunenko, Semisimple restrictions
from ${\rm GL}(n)$ to ${\rm GL}(n-1)$,
{\em J. reine angew. Math.}, {\bf 500} (1998), 83--112.

\bibitem{J.A.Green1}
J.A. Green, Polynomial representations of $GL_n(K)$, Lecture Notes in
  Mathematics, v. 830, {\em Springer-Verlag}, Berlin, Heidelberg, New York,
  1980.

\bibitem{James1}
G.D. James, The representation theory of the symmetric groups, {\em Lecture
  Notes in Mathematics 682}, Springer, Berlin, 1978.

\bibitem{Jantzen1}
J.C. Jantzen, Representations of algebraic groups, Pure and Applied
  Mathematics, 131, {\em Academic Press, Inc.}, Boston, MA, 1987.

\bibitem{Kleshchev_notecs}
A.S. Kleshchev, Completely splittable representations of symmetric groups, {\em
  J. Algebra}, {\bf 181} (1996), n. 2, 584--592.

\bibitem{Shchigolev11}
V.V. Shchigolev, Certain extensions of completely splittable modules, {\em
  Izv. Math.}, {\bf 68} (2004), n. 4, 833--850.

\end{thebibliography}
\end{document}